\documentclass{conm-p-l}  %amsproc
\usepackage{amssymb,amsmath,array,curves}
\usepackage{xypic}

\numberwithin{equation}{section}

\newtheorem{Lemma}{Lemma}[section]\newcommand{\bel}{\begin{Lemma}}\newcommand{\eel}{\end{Lemma}}
\newtheorem{Proposition}[Lemma]{Proposition}\newcommand{\bprop}{\begin{Proposition}}\newcommand{\eprop}{\end{Proposition}}
\newtheorem{Theorem}[Lemma]{Theorem}\newcommand{\bthe}{\begin{Theorem}}\newcommand{\ethe}{\end{Theorem}}
\newcommand{\bpr}{{\bf proof:~}}\newcommand{\epr}{$\blacksquare$\\}
\newtheorem{Remark}[Lemma]{Remark}\newcommand{\beR}{\begin{Remark}\rm}\newcommand{\eeR}{\end{Remark}}
\newtheorem{Definition}[Lemma]{Definition}\newcommand{\bed}{\begin{Definition}}\newcommand{\eed}{\end{Definition}}
\newtheorem{Example}[Lemma]{Example}\newcommand{\bex}{\begin{Example}\rm}\newcommand{\eex}{\end{Example}}
\newtheorem{Corollary}[Lemma]{Corollary}\newcommand{\bcor}{\begin{Corollary}\rm}\newcommand{\ecor}{\end{Corollary}}
\newtheorem{Fact}[Lemma]{Fact}\newcommand{\bfact}{\begin{Fact}\rm}\newcommand{\efact}{\end{Fact}}

\newcommand{\beq}{\begin{equation}}\newcommand{\eeq}{\end{equation}}
\newcommand{\bem}{\begin{displaymath}}\newcommand{\eem}{\end{displaymath}}
\newcommand{\beqa}{\begin{eqnarray}}\newcommand{\eeqa}{\end{eqnarray}}
\newcommand{\bee}{\begin{enumerate}}\newcommand{\eee}{\end{enumerate}}
\newcommand{\bei}{\begin{itemize}}\newcommand{\eei}{\end{itemize}}
\newcommand{\bet}{\begin{tabular}{cccccccc}}\newcommand{\eet}{\end{tabular}}
\newcommand{\bpm}{\begin{pmatrix}}\newcommand{\epm}{\end{pmatrix}}
\newcommand{\ber}{\begin{array}{l}}\newcommand{\eer}{\end{array}}

\newcommand{\di}{\partial}
\newcommand{\ra}{\!\!\rightarrow\!\!}

\newcommand{\tinyM}{\scriptstyle}
\newcommand{\tinyT}{\scriptsize}\newcommand{\smallT}{\footnotesize}
\newcommand{\tinyA}{\tinyM\text{\tinyT}}

\newcommand{\cO}{{\mathcal{O}}}

\newcommand{\mC}{\mathbb{C}}\newcommand{\mD}{\mathbb{D}}

\newcommand{\mP}{\mathbb{P}}\newcommand{\mS}{\mathbb{S}}
\newcommand{\mPN}{\mP_f^{N_d}}\newcommand{\mR}{\mathbb{R}}

\newcommand{\al}{\alpha}\newcommand{\de}{\delta}

\newcommand{\ep}{\epsilon}
\newcommand{\Si}{\Sigma}

\newcommand{\tC}{\tilde{C}}\newcommand{\tSi}{\widetilde{\Sigma}}

\newcommand{\hPl}{\check\mP_l^2}

\newcommand{\lSi}{{\overline{\Sigma}}}\newcommand{\ltSi}{{\overline{\widetilde\Sigma}}}

\newcommand{\li}{~\\ $\bullet$ }

\newcommand{\gNnd}{generalized Newton-non-degenerate }\newcommand{\Nnd}{Newton-non-degenerate }
\newcommand{\ND}{Newton diagram }\newcommand{\D}{diagram }\newcommand{\sqh}{semi-quasi-homogeneous }
\newcommand{\sx}{{\mS_x}}\newcommand{\sy}{{\mS_y}}\newcommand{\sfi}{{\mS_f}}  \newcommand{\omp}{ordinary multiple point}
\newcommand{\Cxy}{{\overset{\frown}{xy}}}  \newcommand{\lxy}{{\overline{xy}}}
\newcommand{\mesh}[7]
{
\put(#1,#2){\vector(1,0){#6}}  \put(#1,#2){\vector(0,1){#7}}
\setcounter{tempx}{#3}  \addtocounter{tempx}{1}  \setcounter{tempy}{#4}  \addtocounter{tempy}{1}
\multiput(#1,#2)(#5,0){\value{tempx}}{\multiput(-1.5,-0.5)(0,#5){\value{tempy}}{.}}
}

\begin{document}\setcounter{secnumdepth}{6} \setcounter{tocdepth}{1}\newcounter{tempx}\newcounter{tempy}
\title{On the collisions of singular points of complex algebraic plane curves}
\author{Dmitry Kerner}
\address{Max Planck Institut f\"ur Mathematik, Vivatsgasse 7,  Bonn 53111, Germany.}
\curraddr{Department of Mathematics, Ben Gurion University of the Negev, P.O.B. 653, Be'er Sheva 84105, Israel.}
\email{kernerdm@math.bgu.ac.il}
\thanks{The research was constantly supported by Max Planck Institut f\"ur Mathematik Bonn, Germany.}
\thanks{The research was also partially supported by the Skirball postdoctoral fellowship of the Center
of Advanced Studies in Mathematics (Mathematics Department of Ben Gurion University) and by
 the Emmy Noether Research Institute for Mathematics and the Minerva Foundation of Germany
}
\date{}
\begin{abstract}
We study the "generic" degenerations of curves with two singular points when the points merge.
First, the notion of generic degeneration is defined precisely. Then a method to
classify the possible results of generic degenerations is proposed in the case of linear singularity types.
We discuss possible bounds on the singularity invariants of the resulting type in terms of the initial types.
In particular the strict upper bound on the resulting multiplicity is proved and a sufficient condition for $\delta=const$
collision is given.
\end{abstract}

\maketitle

\tableofcontents
\section{The problem}$\empty$\\
\parbox{9cm}
{\subsection{Introduction}%$\empty$\\
Let $C$ be a (complex, plane, projective) curve of (high) degree $d$, with singular points $x,y\in\mP^2$ of
(local embedded topological) types $\sx,\sy$. Degenerate $C$ such that the points $x,y$ merge.
}
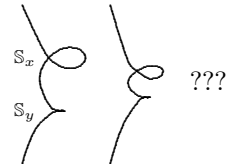
\begin{picture}(0,0)(-20,-15)
\curve(-10,30,-3,15,0,10,4,7,8,6,12,7,14,10,12,13,8,14,4,13,0,10,-3,5,-3,-1,1,-9,4,-10,6,-10)
\curve(6,-10,4,-10,1,-11,-5,-17,-10,-30)
\put(-13,8){$\tinyA \sx$} \put(-13,-12){$\tinyA \sy$}
\end{picture}
\begin{picture}(0,0)(-50,-15)
\curve(-10,30,-3,8,0,5,5,2,10,5,5,7,0,5,-3,2,-3,-1,0,-4,3,-5,5,-5)\curve(5,-5,3,-5,1,-6,-5,-15,-10,-30)
\put(20,-2){???}
\end{picture}
\\\\
We call
this process the {\it collision} of $\sx$ and $\sy$.
What can be said about the {\it resulting singularity type of their generic collision}?

While the collision phenomenon is most natural, it seems to be complicated and not much studied.
In this paper we formulate some specific questions and give some preliminary results.
(For definitions and notations cf. \S\ref{SecAuxilliaryNotionsNotations}.)

To formulate the question precisely, let $\mP H^0(\cO_{\mP^2}(d))=\mPN$ be the parameter space of
plane curves of degree $d\gg0$ (the complete linear system). Here $N_d={d+2\choose{2}}-1$ (the number of monomials of
degree $d$ in 3 variables, minus one).
Consider the subset $\Si_{\sx\sy}\subset\mPN$, the stratum of curves with 2 prescribed singularities.
In this paper we always assume the degree $d$ to be high enough (for the given types $\sx\sy$).
Then the stratum $\Si_{\sx\sy}$ is known to be irreducible, quasi-projective algebraic variety
of expected (co-)dimension.
The boundary of the topological closure ($\lSi_{\sx\sy}\!\!\setminus\!\!\Si_{\sx\sy}$) consists of points
corresponding to curves with higher singularities. In particular, we consider its part
(denoted by $\lSi_{\sx\sy}|_{x=y}$) consisting of all the possible results of collisions of $\sx\sy$.
(It is considered as a reduced subvariety, with the multiplicities omitted.)
\bed
We say that the type $\sfi$ is a result of the collision of types $\sx,\sy$ (and write $\sx+\sy\ra\sfi$) if
a representative of $\sfi$ (i.e. a curve with the singular point of type $\sfi$) belongs to the boundary:
$\lSi_{\sx\sy}|_{x=y}$
\eed
\beR
Note that we do not demand that the {\it whole} stratum $\Si_\sfi$ lie in the boundary $\lSi_{\sx\sy}|_{x=y}$.
In fact examples are known \cite{Pham} (cf. also \cite{DamonGalligo93}) when the adjacency depends on moduli.
So, it is not clear when the collision $\sx+\sy\ra\sfi$ implies the inclusion $\Si_\sfi\subset\lSi_{\sx\sy}|_{x=y}$.
\eeR
The notion "generic" is problematic. The boundary $\lSi_{\sx\sy}|_{x=y}$ is usually reducible, with components
of different dimensions (all of which might be important in applications). One often has to
consider collisions with additional conditions. Say, the tangents $l^x_i$ to
(some of) the branches of $\sx$ (do not) coincide with (some of) those $l^y_j$ of $\sy$. Or, they
(do not) coincide with the limiting tangent line $l=\lxy$ to the curve $\Cxy$, along which the points collide.
In such cases one might be forced
to consider a subvariety of an irreducible component of $\lSi_{\sx\sy}|_{x=y}$.

Therefore, we accept the following definition. For a given singularity type $\mS$, consider the classifying
space of the parameters of the singular germ (e.g. the singular point, the lines of the tangent cone,
with their multiplicities: $T_C=(l_1^{p_1}...l_k^{p_k})$). To a curve with two singular points $\sx\sy$
we assign also the line $l$ through the two points. All this defines a lifting of the initial stratum to a
bigger ambient space:
\beq\ber
\tSi_{\sx\sy}\!\!\!=\!\!\!\Bigg\{\!\!
\Big(\!\!\!\ber(x,\{l^x_i\}..)\\(y,\{l^y_j\}..)\eer\!\!\ber l,C,~x\ne y\\l=\lxy\eer\!\!\!\Big)\Big|~~
C\text{ has}\!\!\ber \text{ $\sx$ at $x$, with $T=((l^x_1)^{p_1}...(l^x_{k_x})^{p_{k_x}})$, ....}\\
\text{ $\sy$ at $y$, with $T=((l^y_1)^{p_1}...(l^y_{k_y})^{p_{k_y}})$, ....}\eer\!\!\!\Bigg\}
\\
\tSi_{\sx\sy}\!\!\!\subset  Aux_x\times Aux_y\times\hPl\times\mPN
\eer\eeq
here $Aux_i$ are the classifying spaces (the notation is for auxiliary), $\hPl$ is the space of lines
in the plane (a line is defined by a one-form). The simplest example is the {\it minimal lifting}
\beq\label{EqMinimalLifting}
\tSi_{\sx\sy}(x,y):=\Big\{\ber(x,y,l,C)\\x\ne y,~l=\lxy\eer\Big|C\text{ has $\sx$ at $x$ and $\sy$ at $y$}
\Big\}\subset\mP^2_x\times\mP^2_y\times\hPl\times\mPN
\eeq
\\\\
Collisions are particular types of degenerations (or deformations). A factorization of the degeneration $\mS_1\ra\mS_3$
is the sequence of degenerations: $\mS_1\ra\mS_2\ra\mS_3$ (i.e. the adjacency of the strata
$\lSi_{\mS_1}\supsetneq\lSi_{\mS_2}\supsetneq\lSi_{\mS_3}$).
\bed The collision $\sx+\sy\ra \sfi$ is called prime if it cannot be factorized (non-trivially).
The collision is called primitive (relatively to a specified lifting)
if the stratum $\tSi_{\sfi}$ contains one of the irreducible components of $\ltSi_{\sx\sy}|_{x=y}$.
\eed
\bex The cases below are well known and can be checked e.g. by the methods of \S \ref{SecAlgorithmCollision}.
\li In the case $A_1+A_1$ we do not fix the tangent lines (as there is no preferred choice). So the lifting
is minimal (eq. \ref{EqMinimalLifting}) and the result is unique: $\ltSi_{A_{1x}A_{1y}}(x,y)|_{x=y}=\ltSi_{A_{3}}$.
The collision is both prime and primitive.
\li In the case $A_2+A_1$ consider first the minimal lifting $\tSi_{A_{2x}A_{1y}}(x,y)$. Then there is only
one primitive (and prime) collision: $\ltSi_{A_{2x}A_{1y}}(x,y)|_{x=y}=\ltSi_{A_4}$.
\\
Now take into account the tangent line of $A_2$  (denoted by $l_x$) and consider the lifting
$\ltSi_{A_{2x}A_{1y}}(x,l_x,y)$. Now, two primitive collisions are possible: $A_2+A_1\ra A_4$ (with $l_x=l$) and
$A_2+A_1\ra D_5$ (with $l_x\neq l$). (The second collision is not prime.)
Indeed the boundary $\ltSi_{A_{2x}A_{1y}}(x,l_x,y)|_{x=y}$  consists of two components: $\ltSi_{A_4}(x,l=l_x)$
and $\widetilde{\ltSi_{D_5}}(x,l,l_x)$. The later is the $\mP^1$ fibration over $\ltSi_{D_5}(x,l)$ with the fiber:
all the lines $l_x$ passing through $x$. In particular the two components are {\it of the same dimension}.

Naively, the first case ($l_x=l$) could be thought of as
the boundary of the second ($l_x\neq l$), but for the minimal lifting the actual situation is
converse (since $\Si_{D_5}\subset\lSi_{A_4}$).
\eex
We have an immediate property of primitivity:
\\
\parbox{10cm}
{\bprop Suppose the lifting is fixed.
 Every non-primitive collision can be factorized through a primitive one (cf. the diagram). Here
$\sfi'\cdots\!\!>\sfi$ is a degeneration (i.e. $\ltSi_{\mS_f}\subsetneq\ltSi_{\mS_f'}$).
\eprop
}
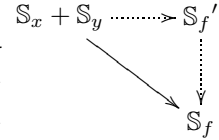
{\xy \POS-(0,-10)\xymatrix{\sx+\sy\ar[dr]\ar@{.>}[r] &\sfi'\ar@{.>}[d]\\&\sfi}\endxy}
\\
The factorization is non-unique, e.g. as one sees from the example above the collision $A_2+A_1\ra D_6$ can be factorized
either through $A_4$ or through $D_5$.
\\
From the last example one sees that the primitivity of the collision depends on the type of lifting.
Thus we {\it fix the choice of lifting for the rest of this paper}.
In the tangent cone of the singularity $T_C=(l_1^{p_1}...l_k^{p_k})$, consider the lines appearing with the
multiplicity 1. They correspond to smooth branches, not tangent to any other branch of the singularity.
We call such branches {\it free}. Call the tangents to the non-free branches: the {\it non-free tangents}.
Assign to the singularity the non-free tangents:
\beq\label{EqLiftingSecondAccepted}\ber
\ltSi_{\sx\sy}\!\!\!:=\!\!\!\overline{\Bigg\{\!\!
\Big(\!\!\!\ber(x,\{l^x_i\})\\(y,\{l^y_j\})\eer\!\!\ber l,C\\x\ne y\eer\!\!\!\Big)\Big\vert
\ber \text{$l^x_i$ are the non-free tangents of $C$ at $x$}\\
\text{$l^y_j$ are the non-free tangents of $C$ at $y$}\eer l=\lxy\Bigg\}}
\\
\ltSi_{\sx\sy}\subset  Aux_x\times Aux_y\times\hPl\times\mPN
\eer\eeq
For \omp s (all the branches are free) this coincides with the minimal lifting.
\beR
\li To specify a collision one should give (at least) the {\it collision data}. It is a list,
specifying the lines among $l,l^x_i,l^y_j$ that merge.
The simplest case is: the limits of all the lines are distinct.
Note that this (seemingly generic) assumption can be often {\it non-generic} (e.g. for the collision $A_{k>1}+A_1$
in the minimal lifting case).
\li We work mostly with {\it linear} singularity types (cf. definition \ref{DefLinearSingularity}).
Typical examples of linear singularities are:
$x^p_1+x^q_2,~~p\le q\le 2p$, $A_{k\le3}$, $D_{k\le6}$, $E_{k\le8}$ etc. Every linear singularity type is necessarily \gNnd
(cf. definition \ref{DefNewtonNonDegenerate}), in particular it has at most two non-free tangents. Even if (some of) the types $\sx\sy$
are non-linear, one can formulate the problem of collisions inside the {\it linear substratum}
$\Si^{(l)}_{\sx\sy}\subset\Si_{\sx\sy}$ (cf. \S\ref{SecSingularityTypes}).
In this case our method provides a complete solution for {\it any} \gNnd singularity type.
\li We always order the types $\sx\sy$ such that $mult_\sx\ge mult_\sy$. In course of collision we always assume $x$ to be fixed.
\eeR
\subsubsection{Acknowledgements}$\empty$\\
Many thanks to G.-M.Greuel, P.Milman and E.Shustin for numerous important discussions.

I am very grateful to the anonymous referee, whose comments helped to improve the text and remove various ambiguities.

The work was done during my stay in Max Planck Institut f\"ur Mathematik, Bonn. I appreciate very much
the excellent working conditions.
\subsection{The specific questions and some partial results}\label{SecSpecificQuestionsPartialResults}
$\empty$\\\subsubsection{A method to classify the results of collision.}$\empty$\\
We propose a method (cf. \S \ref{SecAlgorithmCollision})
to check explicitly the
possible results of a collision, when $\sx$ is \gNnd and $\sy$ is linear.
First we write down the defining equations of the lifted stratum $\tSi_{\sx\sy}$ (outside the diagonal $x=y$).
Then specialize the obtained ideal to the diagonal $x=y$, thus
describing the ideal of the stratum $\ltSi_{\sx\sy}|_{x=y}$. The specialization (the flat limit) is done e.g. by the
usual technic of Gr\"obner basis.
The final step is to recognize the singularity type $\sfi$, from the defining ideal of the stratum $\tSi_\sfi$.

Using the method we discuss in some details the case: $\sy$ is an \omp (\S \ref{SecAlgorithmCollision}).
In particular in the table at the end of paper we list all the possible collision results for the cases:
\li $\sx$ is an \omp (i.e. all its branches are free)
\li one branch of $\sx$ is the ordinary cusp ($x^p_1+x^{p+1}_2$), all others are free (i.e. smooth and non-tangent).
\\\\
There is also a geometric method for some collisions, but these seem to be very special (considered shortly in \S
\ref{SecGeometricApproach}).
$\empty$\\\subsubsection{When does the collision commute with degeneration/deformation?}$\empty$\\
\parbox{6cm}
{Namely, when the diagram on the right commutes? Here the degenerations (deformations) in both rows must
be of course "of the same nature" though applied to the
different types (e.g. \mbox{$A_k\ra A_{k+1}$}, \mbox{$D_k\ra D_{k+1}$}).
}
{\xy \POS-(0,-10)\xymatrix{
\sx^{def}+\sy\ar[d] & \ar[l]_{deform}\sx+\sy\ar[d]\ar[r]^{degen} & \sx^{deg}+\sy\ar[d]
\\
\sfi^{def} & \ar[l]_{deform}\sfi\ar[r]^{degen} & \sfi^{deg}
}\endxy}
\\
We do not know neither how to formulate this question precisely, nor how to answer it.
Nevertheless the idea itself leads to a useful semi-continuity criterion (\S \ref{SecSemiContinuityPrinciple}):
\bprop\label{ThmIntroSemicontinuity}
Let $inv$ be an invariant of the singularity type, upper semi-continuous i.e. non-increasing under the deformations
(i.e. for any small deformation $C_t$ of a curve $C_0$ one has $inv(C_0)\ge inv(C_t)$).
\li Suppose there exists a collision $\sx+\sy\ra\sfi'$ such that $inv(\sfi')$ is the minimal possible for the given case
(e.g. $mult(\sfi')=max(mult(\sx),mult(\sx))$ or $\de_{\sfi'}=\de_\sx+\de_\sy$ or $\mu_{\sfi'}=\mu_\sx+\mu_\sy+1$ etc.).
Then there exists a primitive and prime collision $\sx+\sy\ra\sfi$ such that $inv(\sfi)=inv(\sfi')$.
\li Let $\sx\ra\sx^{deg}$ be a degeneration and $\sx^{deg}+\sy\ra\sfi^{deg}$ a primitive collision. Then there exists
a primitive and prime collision $\sx+\sy\ra\sfi$ and a degeneration, such that the diagram commutes.
In particular, \mbox{$inv(\sfi^{deg})\ge inv(\sfi)$}.
\li Under the assumptions above, suppose for every primitive collision \mbox{$\sx^{deg}+\sy\ra\sfi^{deg}$} the
bound \mbox{$inv(\sfi^{deg})\le a$} is satisfied.
Then for any primitive collision $\sx+\sy\ra\sfi$ one has: $inv(\sfi)\le a$.
\eprop
$\empty$\\\subsubsection{When the collision can be done "branch-wise"?}$\empty$\\
Given the decomposition of both germs $\sx,\sy$ into branches, suppose we know the results of collisions of the branches.
What can be said about the total collision? We give one result in this direction (the proof is
in \S\ref{SecBranchwiseCollisions}):
\bprop\label{ThmBranchwiseCollision} Let a germ $(C,x)$ be the union of two germs: $(C_1,x)\cup(C_2,x)$
(each can be further reducible).
Suppose the germs have no common tangents, i.e. for the tangent cones: $T_{(C_1,x)}\cap T_{(C_2,x)}=\{0\}$.
Denote this by
$\sx=\sx^{\!\!\!1}\cup\sx^{\!\!\!2}$.
\li If there exists a collision $\sx^{\!\!\!2}+\sy\ra\sfi$ then there exists a collision
\mbox{$(\sx^{\!\!\!1}\cup\sx^{\!\!\!2})+\sy\ra(\sx^{\!\!\!1}\cup \sfi)$}
\li If the collision $\sx^{\!\!\!2}+\sy\ra\sfi$ is primitive and $mult(\sfi)=mult(\sx^{\!\!\!2})\ge mult(\sy)$ then the collision
\mbox{$(\sx^{\!\!\!1}\cup\sx^{\!\!\!2})+\sy\ra(\sx^{\!\!\!1}\cup\sfi)$} is primitive.
\eprop

$\empty$\\\subsubsection{Some bounds on the invariants of the resulting types}$\empty$\\
We consider questions of two types. Given the singularity types $\sx\sy$ and a singularity invariant $inv$:
\li what is the ultimate bound on $inv_\sfi$ vs $inv_\sx,inv_\sy$ (i.e. the one satisfied in any primitive collision)?
\li what are the collisions in which $inv$ is "almost unchanged"? (Below we
consider specific invariants and give the precise statements).

The simplest invariants is the multiplicity, here we can answer both questions. More generally,
for {\it semi-continuous} invariants (e.g. $\mu,\de$) there is a hope for some definite answers using the semi-continuity criterion as above.

In the case of non semi-continuous invariants (e.g. the number of branches, the order of determinacy)
we have no hope to obtain any reasonable bounds, we only provide some (counter-)examples.

In this paper we restrict the consideration to the simplest invariants as above. An interesting  question is, of course,
to study the behavior of other invariants (e.g. spectrum, vanishing cohomology etc.)
$\empty$\\\paragraph{What are the possible values of the resulting multiplicity?}$\empty$\\
Here we have two results (proved in \S \ref{SecMultiplicityBounds}):
\bprop\label{ThmIntroBoundOnMultiplicity} Given two types $\sx,\sy$ of multiplicities $m_x\ge m_y$
\li There always exists a collision $\sx+\sy\ra\sfi$ with $mult(\sfi)=m_x$
\li Let $r_x,r_y$ be the number of free branches. If \mbox{$r_x+r_y\ge m_y$} then for any primitive collision
$\sx+\sy\ra\sfi$: $m_\sfi=m_x$.
If $r_x+r_y< m_y$ then for any primitive collision $\sx+\sy\ra\sfi$: $mult(\sfi)\le m_x+m_y-r_x-r_y$.
\eprop
$\empty$\\\paragraph{Bounds on Milnor number}$\empty$\\
A trivial lower bound arises from L\^e-Ramanujam theorem \cite{LeRaman76}: \mbox{$\mu_\sfi\ge\mu_{\sx}+\mu_{\sy}+1$}.
Another bound arises from the semi-continuity of the $\kappa$-invariant, $\kappa=\mu+mult-1$. One has:
$\mu_\sfi\ge\mu_\sx+\mu_\sy+(mult_\sx+mult_\sy-mult_\sfi)-1$.
Probably this bound can be improved:
\\{\bf Conjecture:} If $m_\sx\ge m_\sy$ then  $\mu_\sfi\ge\mu_{\sx}+\mu_{\sy}+m_\sy-1$.
\\The idea of a possible proof is to consider the homological cycles in $H_1(C_{\sx\sy})$ that
vanish as $y\to x$ (cf.\S\ref{SecTopologicalApproach}).
\\\\
An interesting question is to classify the collisions in which the lower bound is realized.

Regarding upper bounds we can only give a very ineffective ones, involving multiplicities and orders of determinacy
(e.g. $\mu_\sfi\le o.d._\sx^2+o.d._\sy^2+o.d._y$)
$\empty$\\\paragraph{The $\de$ invariant (genus discrepancy).}$\empty$\\
How to characterize the $\de=const$ collisions? (They seem to be especially simple.)
This question is partially solved in \S\ref{SecDeltaConstCollisions}. Results there inspire the following
\\{\bf Conjecture:} Given the types $\sx\sy$ let $r_\sx,r_\sy$ be the numbers of potentially free branches (cf.
definition \ref{DefPotentiallyFreeBranches}). Then the $\de=const$ collision is possible iff $r_\sx+r_\sy\ge\min(mult_\sx,mult_\sy)$.
$\empty$\\\paragraph{Other invariants.}\label{SecNumberOfBranches}$\empty$\\
It is much more difficult to give any bounds on non semi-continuous  invariants.
\\
{\bf Number of branches.} We can only give two discouraging examples:
\bex$\bullet$
Consider the primitive collision of two {\it uni-branched} germs $(x^p_1+x^{p+1}_2)$, $(x^2_1+x^{3}_2)$,
 with all 3 lines different $l_x\ne l\ne l_y$.
The resulting type is $(x^{p-2}_1+x^{p-2}_2)(x^3_1+x^4_2)$ (with $p-1$ branches).
\li Consider a primitive collision of two points of the same type: $\sx=\sy=$ \mbox{$\Big(r-1$ free} (smooth) branches
and one cuspidal branch of multiplicity $r-1\Big)$.
Assume that the line of collision is tangent to the cuspidal branches of $\sx,\sy$. It can be shown that the resulting type has only two branches
(of the same type with the common tangent line).
\eex
So a possible upper bound on the number of branches should necessarily involve the multiplicities and probably
there does not exists any lower bound.
\\\\{\bf Order of determinacy.} We can only propose the natural
conjecture:\\ \mbox{$o.d.(\sfi)\le o.d.(\sx)+o.d.(\sy)+1$}.
\\A reasonable lower bound of this conjecture seem to be problematic by the following example.

Consider the collision of two $A_k$ points with distinct tangents. So, the order of determinacy is $k+1$.
It can be shown that among the possible results
(the non-primitive collisions) there is a possibility: $A_k+A_k\ra$ ordinary multiple point of multiplicity
$<\frac{k}{2}+2$. Whose order of determinacy is less than $\frac{k}{2}+2$.

A natural lower bound $o.d.(\sfi)\ge max(o.d.(\sx),o.d.(\sy))$ can be given only for the collision
of linear singularities. In this case it is an immediate consequence of the collision algorithm
(cf. \S\ref{SecAlgorithmCollision}).

$\empty$\\\subsubsection{Topological approach}\label{SecTopologicalApproach}$\empty$\\
The curve $C_{\sx,\sy}$ can be thought of as a partial smoothing of $C_\mS$.
Correspondingly one can deform $C_\mS$ and choose the Milnor fibre so that the vanishing cycles of
$C_{\sx,\sy}$ form a subset of vanishing cycles of $C_\mS$ and the homology lattice for $C_{\sx,\sy}$ is
a sublattice of that for $C_\mS$. Which restrictions does this produce? For example, an $ADE$ singularity
$\mS$ can split to a collection of points of types $\mS_i\in ADE$ iff the union of Dynkin
diagrams $D_{\mS_i}$ can be obtained from $D_\mS$ by deletion of some vertices \cite{Ljashko79}
(cf. also \cite[I.2.7]{AGLV}).
This solves completely the problem of $ADE+ADE\to ADE$ collisions (cf. \S \ref{SecExampleADE}).
A similar statement is known also for some other types of singularities \cite{Jaworski94}.
The natural generalization is therefore:
\\Given the initial types $\sx,\sy$  and a type $\mS$, whose Dynkin diagram $\mD_\mS$ (in some basis) contains
$D_{\sx},D_{\sy}$ (separated by at least one vertex). Is the collision \mbox{$\sx+\sy\to\mS$} possible?

We hope to consider this question in the future.

\section{Auxiliary notions and notations}\label{SecAuxilliaryNotionsNotations}
When considering the local questions, we work in the local coordinates $(x_1,x_2)$ around the point. Working
with the strata we use the homogeneous coordinates $(x_0,x_1,x_2)$ on $\mP^2$.
A (projective)
line through the point $x\in\mP^2$ is defined by a 1-form $l$ (so that $l\in\hPl,~l(x)=0)$.
\\
\\
We often work with symmetric $p-$forms $\Omega^{p}\!\!\in\!\! S^p(\breve{V}_3)$ (here $\breve{V}_3$ is a 3-dimensional vector
space of linear forms).
Thinking of the form as of a symmetric tensor with $p$ indices ($\Omega^{(p)}_{i_1,\dots,i_p}$), we write
$\Omega^{(p)}(\underbrace{x,\dots,x}_{k})$ as a shorthand for the tensor, multiplied $k$ times by the point $x\in\mP_x^2$
\beq
\Omega^{(p)}(\underbrace{x,\dots,x}_{k}):=\sum_{0\le i_1,\dots,i_k\le2}\Omega^{(p)}_{i_1,\dots,i_p}x_{i_1}\dots x_{i_k}
\eeq
So, for example, the expression $\Omega^{(p)}(x)$ is a $(p-1)-$form. Unless stated otherwise, we assume the
symmetric form $\Omega^{(p)}$ to be generic (in particular non-degenerate, i.e. the corresponding curve
$\{\Omega^{(p)}(\underbrace{x,\dots,x}_{p})=0\}\subset\mP^2$ is smooth).

Symmetric forms typically occur as tensors of derivatives of order $p$, e.g. $f^{(p)}$. Sometimes, to emphasize
the point at which the derivatives are calculated we assign it. So, e.g.
$f|_x^{(p)}(\underbrace{y,\dots,y}_k)$ means: the tensor of derivatives of
order $p$, calculated at the point $x$, and contracted $k$ times with $y$.
$\empty$\\\subsection{On the singularity types}\label{SecSingularityTypes}$\empty$\\
\bed\cite{GLSbook}
Let $(C_x,x)\subset(\mC^2_x,x)$ and $(C_y,y)\subset(\mC^2_y,y)$ be two germs of isolated curve singularities.
They are  topologically equivalent if there exist a homeomorphism
 $(\mC^2_x,x)\mapsto(\mC^2_y,y)$ mapping $(C_x,x)$ to $(C_y,y)$. The corresponding equivalence class is
 called the (embedded topological) singularity type. The variety of points (in the parameter space
$\mPN$), corresponding to curves with singularity of a given  (topological) type $\mS$ is called
the {\bf equisingular stratum} $\Si_\mS$
\eed
The topological type can be specified by a (simple, polynomial) representative of the type:
the {\bf normal} form. Several simplest types are
(all the notations are from \cite{AGLV}, we ignore moduli of the analytic classification):
\beq\scriptstyle\ber
A_k:x^2_2+x^{k+1}_1,~~D_k:x^2_2x_1+x^{k-1}_1,~~E_{6k}:x^3_2+x^{3k+1}_1,~~E_{6k+1}:x^3_2+x_2x^{2k+1}_1,\\
E_{6k+2}:x^3_2+x^{3k+2}_1,~~~J_{k\ge1,i\ge0}:x^3_2+x^2_2x^k_1+x^{3k+i}_1,~~Z_{6k-1}:x^3_2x_1+x^{3k-1}_1,\\
Z_{6k}:x^3_2x_1+x_2x^{2k}_1,~~Z_{6k+1}:x^3_2x_1+x^{3k}_1,~~X_{k\ge1,i\ge0}:x^4_2+x^3_2x^k_1+x^2_2x^{2k}_1+x^{4k+i}_1,\\
W_{12k}:x^4_2+x^{4k+1}_1,~~W_{12k+1}:x^4_2+x_2x^{3k+1}_1
\eer\eeq
Using the normal form $f=\sum a_{\bf I}{\bf x}^{\bf I}$ one can draw
the {\it Newton diagram} of the singularity. Namely, one marks the points ${\bf I}$ corresponding to non-vanishing
monomials in $f$, and takes the convex hull of the sets ${\bf I}+\mR_+^2$. The envelope of the convex
hull (the chain of segment-faces) is the Newton diagram.
\bed\cite{GLSbook}\label{DefNewtonNonDegenerate}
\li The singular germ is called \Nnd with respect to its diagram if the truncation of
its polynomial to every face of the diagram is non-degenerate (i.e. the truncated polynomial has no
singular points in the torus $(\mC^*)^2$).
\li The germ is called \gNnd if it can be brought to a Newton-non-degenerate form by a
locally analytic transformation.
\li The singular type is called \Nnd if it has a (generalized) Newton-non-degenerate representative.
\eed
For \Nnd types the normal form is always chosen to be \Nnd.
So, the \Nnd type $\mS$ can be specified by giving the Newton diagram of its normal form $\mD_\mS$.

Newton-non-degeneracy implies strong restrictions on the tangent cone:
\bprop\label{ThmNewtonNonDegenerateTangentCone} Let $T_C=\{(l_1,p_1)...(l_k,p_k)\}$ be the tangent cone of the germ $C=\cup C_j$ (here all the tangents $l_i$
are different, $p_i$ are the multiplicities, so that $\sum_i p_i=mult(C)$). If the germ is \gNnd then $p_i>1$
for at most two tangents $l_i$.
\eprop
So, for a \gNnd germ there are at most two distinguished tangents. We always orient the coordinate axes along
these tangents.

As we consider the topological types, one could expect that to bring a germ to the Newton diagram of the
normal form, one needs local homeomorphisms. However (for curves) the locally analytic transformation always suffice.
In this paper we restrict consideration further to the types for which only linear transformations suffice.
\bed\cite{Ker06}\label{DefLinearSingularity}
A (\gNnd) singular germ is called linear if it can be brought to the Newton diagram of its type by projective
transformations only (or linear transformations in the local coordinate system centered at the
singular point).
A linear stratum is the equisingular stratum, whose open dense part consists of linear germs.
The topological type is called linear if the corresponding stratum is linear.
\eed
The linear types happen to be abundant due to the following observation
\bprop\cite[section 3.1]{Ker06}
The \Nnd topological type is linear iff every segment of the Newton diagram
 has the slope bounded in the segment $[\frac{1}{2},2]$.
\eprop
\bex
The simplest class of examples of linear singularities is defined by the series: $f=x^p+y^q,~~p\leq q\leq2p$.
In general, for a given series only for a few types of singularities the strata can be linear.
In the low modality cases the linear types are:
\li{Simple singularities (no moduli)}: $A_{1\le k\le3},~~D_{4\le k\le6},~~E_{6\le k\le8}$
\li{Unimodal singularities}: $X_9(=X_{1,0}),~~J_{10}(=J_{2,0}),~~Z_{11\le k\le13},~~W_{12\le k\le 13}$
\li{Bimodal}: $Z_{1,0},~~W_{1,0},~~W_{1,1},~~W_{17},~~W_{18}$
\eex
Most singularity types are nonlinear. For example, if a curve has an $A_4$ point, the best we can do by projective
transformations is to bring it to the Newton diagram of $A_3$ $a_{0,2}x^2_2+a_{2,1}x_2x^2_1+a_{4,0}x^4_1$.

This quasi-homogeneous form is degenerated ($a^2_{2,1}=4a_{0,2}a_{4,0}$) and by quadratic (nonlinear!) change
 of coordinates the normal form of $A_4$ is achieved.

Even if a type $\mS$ is non-linear, one can consider the {\it linear substratum}: $\Si^{(l)}_\mS\subset\Si_\mS$ consisting of
points corresponding to all the curves that can be brought to the specified Newton diagram $\mD_\mS$ by linear transformations only.
So, for a linear type $\Si^{(l)}_\mS\equiv\Si_\mS$. Such linear substrata strata are often important (they possess especially nice properties).
\\\\
By the finite determinacy theorem the topological type of the germ is fixed by a finite jet of the defining series. Namely, for every
type $\mS$, there exists $k$ such that for all bigger $n\ge k$: $jet_n(f_1)$ has type $\mS$ iff $f_1$ has type $\mS$.
The minimal such $k$
is called: {\it the order of determinacy}. E.g. $o.d.(A_k)=k+1$, $o.d.(D_k)=k-1$.
The classical theorem is \cite[\S I.2.2]{GLSbook}: if $m^{k+1}\subset m^2Jac(f)$
then $o.d.(f)\le k$.
\section{The results}%$\empty$\\
\subsection{Explicit calculation of collisions}\label{SecAlgorithmCollision}%$\empty$\\
\subsubsection{The idea}$\empty$\\
One way of treating the problem could be to consider explicit equations of the stratum $\ltSi_{\sx\sy}$ and
then to restrict them to the diagonal $x=y$. But it is difficult to write down the complete set of
the generators of the ideal $I(\ltSi_{\sx\sy})$. Instead, we start from the ideals $I(\ltSi_{\sx})$, $I(\ltSi_{\sy})$
of the coordinate ring $K[Aux_x\times Aux_y\times\hPl\times\mPN]$.
Their sum $I(\ltSi_{\sx})\oplus I(\ltSi_{\sy})$ defines the stratum $\ltSi_{\sx\sy}$ outside the diagonal.
Over the diagonal the sum does not define the stratum (since the intersection $\ltSi_{\sx}\cap \ltSi_{\sy}$
has residual components of excess dimension).

One way to continue is to take the topological closure: $\overline{\ltSi_{\sx}\underset{x\ne y}{\cap}\ltSi_{\sy}}$.
From the calculational point of view we should take the flat limit of $I(\ltSi_{\sx})\oplus I(\ltSi_{\sy})$
as $y$ approaches $x$.

More formally, we use the standard fact:
\bprop
The flat limit of $I(\ltSi_{\sx})\oplus I(\ltSi_{\sy})$ as $y\ra x$ gives the defining ideal $I(\ltSi_{\sx\sy}|_{x=y})$.
\eprop
To take the flat limit, one should preserve all the inter-relations (syzygies).
This is done e.g. by finding the Gr\"obner basis \cite[section 2]{StevensBooks}.

Thus the problem is reduced (at least theoretically) to the study of ideals $I(\ltSi_{\sx})$, $I(\ltSi_{\sy})$.
For many singularity types the generators of the ideals are known \cite{Ker06} and can be written in a simple form.
These types include the {\it linear} singularities (cf. the definition \ref{DefLinearSingularity}). Examples of such types are
$A_{k\le3}$, $D_{k\le6}$, $E_{k\le8}$, $x^p_1+x^q_2,$ $p\le q\le 2p$...

In fact we attack a more general case: when the type $\sy$ is linear and $\sx$ is \gNnd.
Start from a \gNnd type $\sx$, bring the corresponding germ to a \Nnd form by a locally analytic transformation.
 Since the
result of collision is invariant under the locally analytic transformations of $\mC^2$,
can assume that the germ $\sx$ is brought to its \ND by linear transformations. Consider the corresponding subvariety
$\ltSi^l_{\sx}\subset\ltSi_{\sx}$ consisting of those
germs that can be brought to their \ND by {\it linear} transformations. (In \cite{Ker06} such a subvariety was called
the linear substratum.) If the type $\sx$ is itself linear then of course $\ltSi^l_{\sx}\equiv\ltSi_{\sx}$.
So, all the collision results can be (and will be) obtained from $I(\ltSi^{(l)}_\sx)$, $I(\ltSi_\sy)$.
$\empty$\\\subsubsection{How to simplify the collision trajectory}\label{SecCollisionTrajectory}$\empty$\\
We always keep the point $x$ and at least one of the non-free tangents to $(C,x)$ fixed. In general $y$ approaches $x$
along a (smooth) curve $\Cxy:$ given by $y=x+\sum \ep^iv_i$. To simplify the problem, one would like to rectify the curve
 into the line $l=\lxy$
(by a locally analytic transformation preserving the tangents). But our method places severe restrictions on the possible transformations.
Recall that we assume $\sy$ to be a linear type, while $\sx$ is \gNnd. To be able to write the defining conditions, the germ $(C,x)$ is assumed to possess
the \ND of the type $\sx$.
\li If $\sx$ is linear then all the transformations preserving the tangents are allowed (i.e. $x_i\ra x_i+\phi_i,~~\phi\in m^2$).
In particular, the collision can always be assumed to happen along a line.
\li If $\sx$ is not linear (but \gNnd), then only the transformations preserving the diagram are allowed.
So, if the tangent to $\Cxy$ is
distinct from all the non-free tangents of $\sx$, then the curve $\Cxy$ can be rectified to the line $\lxy$. Otherwise, one can
only get an upper bound on the degree of the curve $\Cxy$.
$\empty$\\\subsubsection{The algorithm}$\empty$\\
The input, initial data, consists of the two strata $\lSi^l_{\sx},\lSi_{\sy}$, with known generators of their ideals:
\beq
I(\ltSi_{\sx})=<\{h_i(x)\}_i>,~~~~~~I(\ltSi_{\sy})=<\{g_i(y)\}_j>
\eeq
Here the points $x,y$ are assigned to emphasize the dependence. (Of course, the generators depend on other
 parameters of the singularity also.) Fix the collision data of the types $\sx,\sy$:  $l^x_i,l^y_j,l$.
$\empty$\\\paragraph{Preparation of the series.}$\empty$\\
Expand $y=x+\sum_i \ep^i v_i$. Here $\ep$ is an infinitesimal parameter, while the vectors $v_i$ define the trajectory
of collision. The collision in general happens along a (smooth) curve and higher order expansion parameters
of the curve can be important (e.g. this is the case in $A_{k\ge4}+A_1$ collision).
Expand, all the generators $g_j(y)$ into power series of $\ep$, i.e. $g_j(y)=g_j(x)+\ep()+...$.
Restrict to $\ltSi_\sx$, i.e. take into account the equations of $\sx$. Depending on the collision data, some additional terms in the
series $g(y)$ can vanish.
$\empty$\\\paragraph{Taking the flat limit.}$\empty$\\
 Given the ideal generated by polynomials $f_i(x)$ and series $g_j(y)=g_j(x)+\ep...$ check
all the relations among $\{f_i(x),g_j(x)\}$, i.e. calculate the syzygies.  For every such a relation
$\sum r_i(x)f_i(x)+\sum R_j(x)g_j(x)=0$ one gets a new series $\sum r_i(x)f_i(x)+\sum R_j(x)g_j(y)$. By construction
this series has the common factor: a power of $\ep$.

We work outside the diagonal
(in the ring $K[[Aux_x\times Aux_y\times\hPl\times\mPN,\ep,\ep^{-1}]]$). Therefore, each time
one gets a series with a common factor of $\ep$, divide by $\ep$.

Add all the new series to the initial ideal and check for the new relations (syzygies).
By the general theory, after a finite number of steps the procedure terminates: the standard (Gr\"obner )
basis is constructed.
Now take the limit $\ep\ra0$, omitting all the higher order terms. The obtained system is the system
of generators of the ideal $I(\ltSi_{\sx\sy}|_{x=y})$.

The variety $\ltSi_f$ might be reducible (or non-reduced), in this case take a reduced irreducible component.

The process depends in general on the (non-)coincidence of various tangents to the branches, the collision
line $\bar{xy}$ (i.e. the tangent to the collision curve), the conic osculating to the collision curve etc.

Note that the initial system of generators $f_i(x),g_j(y)$ has a lot of structure (cf. the example
\S \ref{SecCollisionWithMultiplePoint}),
various equations are combined into some symmetric forms. Preserving this structure helps to recognize the
resulting types.

$\empty$\\\paragraph{Recognition of the final singularity type $\sfi$.}\label{SecRecognitionOfTheFinalType}$\empty$\\
As the result of the flat limit (above) we obtain the defining ideal of $\ltSi_f$.
This gives the defining equations of $\ltSi_f$, written in terms of the coordinate $x$, the function $f$
(and its derivatives), the tangents $l^x_i, l^y_j$ to $\sx,\sy$ and parameters of the collision trajectory
$y=x+\sum \ep^iv_i$. If some of the tangents
$l^x_i, l^y_j$ coincide, then we should also consider the way they approach: $l^y_j=l^x_j+\sum \ep^iw_i$.
\bprop
The resulting system of equations is linear in $f$ (and its derivatives).
\eprop
\bpr Note, that all the initial equations $f_i(x),g_j(y)$ are linear in $f$ (since we work with linear (sub)strata)
and are homogeneous in other variables. So, if $f_1,f_2\in\mPN$ satisfy the system then any linear
combination $\al_1 f_1+\al_2f_2$ satisfies it (for other variables fixed, and $x\ne y$). Thus there can
be no relation among the equation involving the function (or its derivatives), except for a trivial
one (Koszul). Since this would produce a non-trivial equation non-linear in $f$.
\epr
\li The simplest case is when the initial system involves only $l^x_i, l^y_j,v,f$ (e.g. both $\sx\sy$ are linear).
Then, as follows from the proposition the resulting stratum is {\it linear}.
Thus the singular type is easy to recognize (can write down a particular simple representative, to draw the \ND etc.).
\li When parameters of the expansions $y=x+\sum \ep^iv_i$, $l^y_j=l^x_j+\sum \ep^iw_i$ appear explicitly in
the equations, the situation is more complicated (i.e. the resulting stratum might be non-linear).
One possible way is to fix some specific values of the parameters and find a specific (generic) solution.
By the proposition above the equations are still linear in $f$, so there is no problem finding an explicit solution.
From this solution one can construct e.g. the resolution tree and thus identify the type.
\beR
Probably in this way one can get some information about the final Newton diagram. Unfortunately
we do not have any result by now.
\eeR

$\empty$\\\subsubsection{Application: a \gNnd singularity $\sx$ and the \omp\\\mbox{$\sy=x^{q+1}_1+x^{q+1}_2$}}\label{SecCollisionWithMultiplePoint}$\empty$\\
Here we assume $mult(\sx)=p+1\ge mult(\sy)=q+1$ and the collision data is generic, i.e. the curve $\Cxy$ is
not tangent to any of the non-free branches of $\sx$. Thus (cf. \S \ref{SecCollisionTrajectory}) the curve $\Cxy$ can be assumed to
be a line: $\Cxy=\lxy=l$.

We should translate the conditions at the point $y$ to conditions at $x$.
Outside the diagonal $x=y$ the stratum is defined by the set of conditions corresponding to $\ltSi_{\sx}$, and by
the condition $f|_y^{(q)}=0$. This is the (symmetric) form of derivatives of order $q$, calculated at the
point $y$ (in projective coordinates). In the neighborhood of $x$ expand $y=x+\sum_i\ep^i v_i$ (here $\ep$ is
small and $v_1$ is the direction along the line $l=\bar{xy}$).
Since we have assumed that the collision happens along a line, in the above expansion we need only the first term:
$y=x+\ep v$.

To take the flat limit, expand $f|_y^{(q)}$ around $x$, we get
$0=f|_y^{(q)}=f|_x^{(q)}+..+\frac{\ep^{p-q}}{(p-q)!}f|_x^{(p)}(v..v)+..$.
First several terms in the expansion vanish, up to the multiplicity of $\sx$. Normalize by the common factor of $\ep$:
\beq\label{EqInitialSystemOrdMultPoin}
\frac{1}{(p-q+1)!}f|_x^{(p+1)}(\underbrace{v..v}_{p+1-q})+\frac{\ep}{(p-q+2)!} f|_x^{(p+2)}(\underbrace{v..v}_{p+2-q})+
\frac{\ep^2}{(p-q+3)!} f|_x^{(p+3)}(\underbrace{v..v}_{p+3-q})+\dots
\eeq
To take the flat limit, we should find all the syzygies between these series and the equations for $\ltSi_{\sx}$.
First we find the "internal" syzygies of the series themselves.
\bprop
The standard basis, obtained by considering all the syzygies of the equation (\ref{EqInitialSystemOrdMultPoin}), is:
\beq\label{EqTriangularSystemOrdMultPoint}%\smallT
\begin{tabular}{@{}>{$}c<{$}@{}>{$}c<{$}@{}>{$}c<{$}@{}>{$}c<{$}@{}>{$}c<{$}@{}>{$}c<{$}@{}>{$}c<{$}@{}>{$}c<{$}@{}>{$}c<{$}@{}>{$}c<{$}}
%\begin{tabular}{@{}>{$}c<{$}@{}>{\smallT$}c<{$}@{}>{\smallT$}c<{$}@{}>{\smallT$}c<{$}@{}>{\smallT$}c<{$}@{}>{\smallT$}c<{$}@{}>{\smallT$}c<{$}@{}>{\smallT$}c<{$}@{}>{\smallT$}c<{$}@{}>{\smallT$}c<{$}}
\\
 f|_x^{(p+1)}(\underbrace{v..v}_{p+1-q}) &
+& \ep f|_x^{(p+2)}(\underbrace{v..v}_{p+2-q})&
+&\ep^2f|_x^{(p+3)}(\underbrace{v..v}_{p+3-q})&
+&\ep^3f|_x^{(p+4)}(\underbrace{v..v}_{p+4-q})&+..
\\
0 &+&  f|_x^{(p+2)}(\underbrace{v..v}_{p+3-q})&
+&\ep f|_x^{(p+3)}(\underbrace{v..v}_{p+4-q})&
+&\ep^2 f|_x^{(p+4)}(\underbrace{v..v}_{p+5-q})&+..
\\
0 &+& 0&+& f|_x^{(p+3)}(\underbrace{v..v}_{p+5-q})
&+& \ep f|_x^{(p+4)}(\underbrace{v..v}_{p+6-q})&+..
\\
..&..&..&..&
\\
0&+&0 &+&0&+&...&+&   f|_x^{(p+q+1)}(\underbrace{v..v}_{p+q+1})&+..
\end{tabular}
\eeq
\eprop
\bpr
The syzygies are obtained as a consequence of the Euler identity for homogeneous polynomial
\mbox{$\sum x_i\di_if=deg(f)f$}.
By successive contraction of the tensor series with $x$ we get the series
\beq\smallT
\begin{tabular}{@{}>{$}c<{$}@{}>{$}c<{$}@{}>{$}c<{$}@{}>{$}c<{$}@{}>{$}c<{$}@{}>{$}c<{$}@{}>{$}c<{$}@{}>{$}c<{$}}
%\begin{tabular}{@{}>{\tinyT$\tinyM}c<{$}@{}>{\tinyT$\tinyM}c<{$}@{}>{\tinyT$\tinyM}c<{$}@{}>{\tinyT$\tinyM}c<{$}@{}>{\tinyT$\tinyM}c<{$}@{}>{\tinyT$\tinyM}c<{$}@{}>{\tinyT$\tinyM}c<{$}@{}>{\tinyT$\tinyM}c<{$}}
\\
 \frac{1}{(p-q+1)!}f|_x^{(p+1)}(\underbrace{v..v}_{p+1-q})&
+& \frac{\ep}{(p-q+2)!} f|_x^{(p+2)}(\underbrace{v..v}_{p+2-q})&
+& \frac{\ep^2}{(p-q+3)!} f|_x^{(p+3)}(\underbrace{v..v}_{p+3-q})&
%+& \frac{\ep^3}{(p-q+4)!} f|_x^{(p+4)}(\underbrace{v..v}_{p+4-q})&+..
\\
\frac{(d-p-2)}{(p-q+2)!}f|_x^{(p+1)}(\underbrace{v..v}_{p+2-q}) &
+& \frac{\ep (d-p-3)}{(p-q+3)!}f|_x^{(p+2)}(\underbrace{v..v}_{p+3-q})&
+& \frac{\ep^2(d-p-4)}{(p-q+4)!} f|_x^{(p+3)}(\underbrace{v..v}_{p+4-q})&
%+& \frac{\ep^3(d-p-5)}{(p-q+5)!} f|_x^{(p+4)}(\underbrace{v..v}_{p+5-q})&+..
\\
%(d-p-2)(d-p-3)f|_x^{(p+1)}(\underbrace{v..v}_{p+3-q}) +& \ep (d-p-3)(d-p-4)f|_x^{(p+2)}(\underbrace{v..v}_{p+4-q}) +& \ep(d-p-4)(d-p-5) f|_x^{(p+3)}(\underbrace{v..v}_{p+5-q})+&..
%\\
..&..&..&..&
\\
\frac{\prod^{q+1}_{i=2}(d-p-i)}{(p+1)!}f|_x^{(p+1)}(\underbrace{v..v}_{p+1})&
+&\frac{\ep\prod^{q+1}_{i=2}(d-p-1-i)}{(p+2)!} f|_x^{(p+2)}(\underbrace{v..v}_{p+2})&
+&\frac{\ep^2\prod^{q+1}_{i=2}(d-p-2-i)}{(p+3)!} f|_x^{(p+3)}(\underbrace{v..v}_{p+3})&
%+& \frac{\ep^3\prod^{q+1}_{i=2}(d-p-3-i)}{(p+4)!} f|_x^{(p+4)}(\underbrace{v..v}_{p+4})&+..
\end{tabular}
\eeq
\\
Here the first row is the initial series, the second is obtained by contraction with $x$ once, the $p+2$'th
row is obtained by contracting ($p+1$) times with $x$.

Apply now the Gaussian elimination, to bring this system to the upper triangular form.
\li Eliminate from the first column all the entries of the rows $2..(p+2)$. For this contract the first
row sufficient number of times with $v$ (fix the numerical coefficient) and subtract.
\li Eliminate from the second column all the entries of the rows $3..(p+2)$.
\li ...
\\Normalize the rows (i.e. divide by the necessary power of $\ep$).

In this way we get the "upper triangular" system of series in eq. (\ref{EqTriangularSystemOrdMultPoint})
(we omit the numerical coefficients).

There are no more "internal" syzygies, i.e. we have obtained the Gr\"obner basis for
the initial system (\ref{EqInitialSystemOrdMultPoin}).
\epr

Now the generators of $I(\ltSi_{\sx})$ should be added and one checks again for the possible syzygies.
As the simplest example consider the case $\sx=x^{p+1}_1+x^{p+1}_2$, $\sy=x^{q+1}_1+x^{q+1}_2$.
\bcor
Let $\sx\sy$ be the ordinary multiple points with multiplicities $(p+1)$ and $(q+1)$ respectively (with $p\ge q$).
For the lifting $\ltSi_{\sx\sy}(x,y)$ there exists only one primitive collision $\sx+\sy\ra\sfi$ with the final
type having the normal form $(x^{p-q}_1+x^{p-q}_2)(x^{q+1}_1+x^{2q+2}_2)$.
\ecor
\bpr
The defining equations of the stratum $\ltSi_\sx(x)$ are: $f|_x^{(p)}=0$
(as there are no non-free branches the lifting is minimal).
Therefore in equation (\ref{EqTriangularSystemOrdMultPoint}) there are no more syzygies, so
 just take the limit $\ep\to0$ (i.e. omit the higher order terms in each row). Finally, we get the defining
system of equations:
\beq{\smallT\tinyM
f|_x^{(p)}=0,~~~f|_x^{(p+1)}(\underbrace{v..v}_{p+1-q})=0,~~~f|_x^{(p+2)}(\underbrace{v..v}_{p+3-q})=0,
~~f|_x^{(p+3)}(\underbrace{v..v}_{p+5-q})=0~~..,f|_x^{(p+q+1)}(\underbrace{v..v}_{p+q+1})=0
}\eeq
\parbox{9cm}{
As was emphasized in \S\ref{SecRecognitionOfTheFinalType}, the system is linear in $f$, so it defines a linear
(sub)stratum $\ltSi^l_{\sfi}$.
We can obtain the \ND of the resulting type by fixing (in projective coordinates) e.g. $x=(0,0,1)$, $v=(0,1,0)$.
Since all the slopes of the \D are lie in the segment $[\frac{1}{2},2]$ we get that the type is linear
and $\ltSi^l_{\sfi}=\ltSi_{\sfi}$.
}
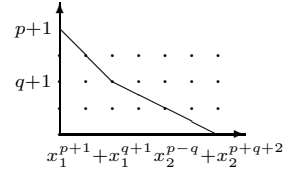
\begin{picture}(0,0)(-20,10)
\mesh{0}{0}{6}{3}{10}{70}{50}
%% Draws the\mesh{x-coordinate of(0,0),y-coordinate of(0,0),number of x-points(0 isn't counted),number of y-points,delta,x-length,y-length}
\put(0,40){\line(1,-1){20}} \put(20,20){\line(2,-1){40}}
\put(-17,38){$\tinyA p+1$}  \put(-17,18){$\tinyA q+1$} %\put(12,-8){$\tinyA p-q$}  \put(45,-8){$\tinyA p+q+2$}
\put(-5,-10){$\tinyA x^{p+1}_1+x^{q+1}_1x^{p-q}_2+x^{p+q+2}_2$}
\end{picture}
\\
\epr
In several simplest cases we have: $A_1+A_1\ra A_3$, $D_4+A_1\ra D_6$, $X_9+A_1\ra X_{1,2}$, $D_4+D_4\ra J_{10}$,
 $X_9+D_4\ra Z_{13}$.
\\\\
\bex{$\sx=x^{p+1}_1+x^{p+2}_2$.} Now the result of collision depends on the (non)coincidence of the line $l=\lxy$
with the tangent line $l_x$ to $\sx$. The lifted stratum $\ltSi_\sx$ is defined by the condition (cf. \cite{Ker06})
$f|_x^{(p+1)}\sim \underbrace{l_x\times..\times l_x}_{p+1}$, this can be written also as
$f|_x^{(p+1)}(v_x)=0$.
\bprop
For $l\ne l_x$ the only (primitive) resulting type is\\ \mbox{$(x^{p+1-q}_1+x^{p+1-q}_2)(x^{q+1}_1+x^{2q+1}_2)$}.
For $l=l_x$
the only (primitive) resulting type is $(x^{p-q}_1+x^{p+1-q}_2)(x^{q+1}_1+x^{2q+2}_2)$.
\eprop
\bpr As $\sx,\sy$ are linear, can assume that the trajectory is a line: $l=\lxy$.
\li $l_x\ne l$. Contract the first row of (\ref{EqTriangularSystemOrdMultPoint}) with $v_x$. The $\ep^0$ term vanish and
the whole series is divided by $\ep$. So, we get: $0=f|_x^{(p+2)}(\underbrace{v..v}_{p+2-q}v_x)+...$. Contract this series with $v$
and subtract from the third row of (\ref{EqTriangularSystemOrdMultPoint})  (contracted with $v_x$). Apply the same procedure,
up to the last row. Direct check shows that there are no more syzygies, so substitute $\ep=0$
and get
\beq\ber\smallT
f|_x^{(p+1)}\sim \underbrace{l_x\times..\times l_x}_{p+1},~~~f|_x^{(p+1)}(\underbrace{v..v}_{p+1-q})=0,~~~
f|_x^{(p+2)}(\underbrace{v..v}_{p+3-q})=0,~~f|_x^{(p+2)}(\underbrace{v..v}_{p+2-q},v_x)=0,\\
f|_x^{(p+3)}(\underbrace{v..v}_{p+5-q})=0,~~f|_x^{(p+3)}(\underbrace{v..v}_{p+4-q}v_x)=0,~~..,
f|_x^{(p+q+1)}(\underbrace{v..v}_{p+q+1})=f|_x^{(p+q+1)}(\underbrace{v..v}_{p+q}v_x)=0
\eer\eeq
which gives (since $v_x\ne v$ and $l_x(v)\ne0$):\\
\beq\smallT
f|_x^{(p+1)}=0,~~~f|_x^{(p+2)}(\underbrace{v..v}_{p+2-q})=0,f|_x^{(p+3)}(\underbrace{v..v}_{p+4-q})=0,\dots,f|_x^{(p+q+1)}(\underbrace{v..v}_{p+q})=0
\eeq
\parbox{8cm}
{From here we get the normal form: \mbox{$(x^{p+1-q}_1+x^{p+1-q}_2)(x^{q+1}_1+x^{2q+1}_2)$}.
}
\begin{picture}(0,0)(-20,20)
\mesh{0}{0}{6}{3}{10}{90}{50}
%% Draws the\mesh{x-coordinate of(0,0),y-coordinate of(0,0),number of x-points(0 isn't counted),number of y-points,delta,x-length,y-length}
\put(0,40){\line(1,-1){20}} \put(20,20){\line(3,-1){60}}   %\put(40,10){\line(4,-1){40}}
\put(-17,38){$\tinyA p+2$}\put(-17,18){$\tinyA q+1$} %   \put(-17,8){$\tinyA q-1$}
%\put(17.5,17.5){$\bullet$} \put(37.5,7.5){$\bullet$}
%\put(8,-8){$\tinyA p+1-q$} \put(28,-15){$\tinyA p+4-q$} \put(70,-8){$\tinyA p+q+2$}
\put(-15,-12){$\tinyA x^{p+2}_1+x^{q+1}_1x^{p+1-q}_2+x^{p+q+2}_2$}
\end{picture}
\\\\\\
\li $l_x=l$. In this case the system should be re-derived, starting from eq. (\ref{EqInitialSystemOrdMultPoin}).
Everything is just shifted ($p\ra p+1$) and we get the equations:
\beq\ber
f|_x^{(p)}=0,~~~f|_x^{(p+1)}(\underbrace{v..v}_{p+1-q})=0,~~~f|_x^{(p+2)}(\underbrace{v..v}_{p+2-q})=0,\\
f|_x^{(p+3)}(\underbrace{v..v}_{p+4-q})=0~~..,f|_x^{(p+q+2)}(\underbrace{v..v}_{p+q+2})=0
\eer\eeq
\parbox{9cm}
{This gives the normal form of the singularity \mbox{$(x^{p-q}_1+x^{p+1-q}_2)(x^{q+1}_1+x^{2q+2}_2)$}
}
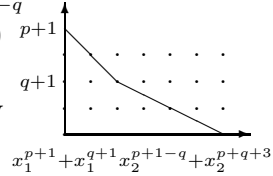
\begin{picture}(0,0)(-20,0)
\mesh{0}{0}{6}{3}{10}{70}{50}
%% Draws the\mesh{x-coordinate of(0,0),y-coordinate of(0,0),number of x-points(0 isn't counted),number of y-points,delta,x-length,y-length}
\put(0,40){\line(1,-1){20}} \put(20,20){\line(2,-1){40}}
\put(-17,38){$\tinyA p+1$} \put(-17,18){$\tinyA q+1$}
%\put(12,-8){$\tinyA p+1-q$}  \put(45,-8){$\tinyA p+q+3$}
\put(-20,-12){$\tinyA x^{p+1}_1+x^{q+1}_1x^{p+1-q}_2+x^{p+q+3}_2$}
\end{picture}
\eex
\epr
$\empty$\\\subsubsection{More general case}$\empty$\\
If the curve $\Cxy$ is tangent to one of the non-free branches of $\sx$, then the system
(\ref{EqTriangularSystemOrdMultPoint})
should be re-derived. When $\sx$ is linear, we can assume that $\Cxy=\lxy=l$, this greatly simplifies the calculations.

If $\sx$ is not an \omp, then to the conditions
of the system (\ref{EqTriangularSystemOrdMultPoint}), one adds the
conditions of $\sx$ and checks for possible additional syzygies.

In some cases there are no new syzygies. For example, let the tangent cone of $\sx$, with
multiplicities be \mbox{$T_{C_x}=\{l^{p_1}_1..l^{p_k}_k\}$}, such
that $\forall i:~~ p_i\le p+1-q$. Consider the primitive collision $\sx+\sy\ra\sfi$ such that the
collision line $l$ is distinct from all the tangents with $p_i>1$. Then the defining ideal of
the resulting stratum is especially simple:
\beq
I(\lSi_{\sfi})=<I(\lSi_{\sx}),\underset{\ep\to0}{\lim}I_\ep>
\eeq
here  $I_\ep$ is the ideal of the equation (\ref{EqTriangularSystemOrdMultPoint}).
$\empty$\\\subsubsection{Geometric approach}\label{SecGeometricApproach}$\empty$\\
\parbox{9cm}
{A natural idea is to trace the collision explicitly by drawing a (real) picture. A useful trick is to blow up the plane at $x$.
}
\begin{picture}(0,0)(-20,-10)
\curve(0,10,10,-5,20,10)\curve(0,-10,10,5,20,-10)\put(-7,-2){$\tinyA A_1$}\put(17,-2){$\tinyA A_1$}\put(30,-5){$\Rrightarrow$}
\end{picture}
\begin{picture}(0,0)(-60,-10)
\curve(0,10,10,-0,20,10)\curve(0,-10,10,0,20,-10)\put(4,4){$\tinyA A_3$}
\end{picture}
\\
If needed one might blowup several times (for example resolve the germ $(C,x)$).
Then make a choice by gluing the branches of ${\sx},\sy$ and collide (i.e. push $\sy$ to the exceptional divisor).
 In this way some parts of the curve are contracted.
\\
\parbox{7cm}
{One gets a curve on the blown up plane,
with a singular point on $E$. Now, blow down (i.e. contract the exceptional divisor).
This gives the resulting germ.
\bex
The collision of two \omp s. Suppose, the multiplicities of $\sx,\sy$ are $p+1,q+1$ such that $p\ge q$.
Blowup at $x$, push $y$
to the exceptional divisor, then blowdown, as in the picture.
\eex
}
\begin{picture}(0,0)(-30,10)
%Left lower part
\qbezier(0,10)(0,0)(0,-10)  \qbezier(-4,10)(0,0)(4,-10) \qbezier(4,10)(0,0)(-4,-10)
 \qbezier(-10,10)(20,-27)(60,10)    \qbezier(-10,-10)(20,27)(60,-10)   \qbezier(-10,0)(0,0)(60,0)
\put(-5,-15){\tinyT p}  \put(45,-10){\tinyT q}  \put(-3,15){$\downarrow$}
%%%Left upper part
\put(0,25){\line(0,1){50}} \put(3,74){\tinyT E} \qbezier(-10,30)(0,30)(10,30) \qbezier(-10,70)(0,70)(10,70)
 \qbezier(-10,40)(40,40)(60,60)   \qbezier(-10,50)(40,50)(60,50)   \qbezier(-10,60)(40,60)(60,40)
  \put(-20,50){\tinyT p}  \put(45,40){\tinyT q}
%%RightUpper part
 \put(70,45){$\rightarrow$}
\put(100,25){\line(0,1){50}} \put(103,74){\tinyT E} \qbezier(90,30)(100,30)(110,30) \qbezier(90,70)(100,70)(110,70)
 \qbezier(90,40)(100,50)(110,60)   \qbezier(90,50)(100,50)(110,50)   \qbezier(90,60)(100,50)(110,40)
%%RightLowerPart
 \put(97,15){$\downarrow$}
\qbezier(90,0)(100,0)(110,0)
 \qbezier(90,-10)(100,10)(110,-10)  \qbezier(90,10)(100,-10)(110,10)
\qbezier(95,13)(100,0)(105,-13)   \qbezier(95,-13)(100,0)(105,13)
\end{picture}
\\
\parbox{8cm}
{More generally, suppose the number of free branches for the type $\sx$ is at least the multiplicity of $\sy$.
Use the same procedure as above, to get the final answer.
}
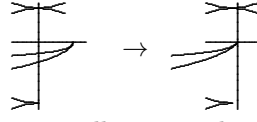
\begin{picture}(0,0)(-20,-10)
\qbezier(0,20)(0,0)(0,-20)  \qbezier(-10,-16)(-5,-18)(0,-18) \qbezier(-10,-20)(-5,-18)(0,-18)
 \qbezier(-10,16)(0,20)(10,16)    \qbezier(-10,20)(0,16)(10,20)

\curve(-10,-5,7,0,13,5)  \curve(-10,0,8,2,13,5)  \curve(-10,5,18,5)
\put(35,0){$\ra$}
\end{picture}
\begin{picture}(0,0)(-90,-10)
\qbezier(2,20)(2,0)(2,-20)  \qbezier(-10,-16)(-5,-18)(0,-18) \qbezier(-10,-20)(-5,-18)(0,-18)
 \qbezier(-10,16)(0,20)(10,16)    \qbezier(-10,20)(0,16)(10,20)

\curve(-23,-5,-6,0,2,5)  \curve(-23,0,-5,2,2,5)  \curve(-10,5,10,5)
\end{picture}
\\
 The restrictions of this approach are evident:
 the primitive collision can be traced for some special types only.
 The curve families which can be simultaneously blown up are usually equi-normalizable, thus
by the classical Teissier theorem (cf. \ref{ThmTeissier}) in such a collision $\de=const$.

In addition, working with real pictures we necessarily loose information.
Here an important fact is that to perform the $\de=const$ collision one can always choose real representatives
of the type (\cite{ACampo75-1,ACampo75-2}, \cite{Gusein-Zade 74-1,Gusein-Zade 74-2})
$\empty$\\\subsubsection{Branch-wise collisions}\label{SecBranchwiseCollisions}$\empty$\\
Here we prove the proposition \ref{ThmBranchwiseCollision}.
\\
\bpr
$\bullet$ The existence of collision can be easily seen e.g. by geometric consideration. Blowing-up the plane at $x$ separates the germs
$(C_1,x)$ and $(C_2,x)$. Thus on the blown up plane can do the collision of the transform of $\sx^{\!\!\!\!2}$ with $\sy$. Now blowdown.
\li Suppose the collision is non-primitive then it can be factorized: $(\sx^{\!\!\!\!1}\cup\sx^{\!\!\!\!2})+\sy\ra\sfi'\stackrel{degen}{\to}(\sx^{\!\!\!\!1}\cup\sfi)$.
Here by the assumptions of the proposition the degeneration should preserve the multiplicity. Therefore the tangent cone $T_{\sx^{\!\!\!\!1}\cup\sfi}$ is
the degeneration limit of $T_{\sfi'}$. So, $\sfi'$ also has (at least) two subsets of branches: $\sfi'=\sfi'_1\cup\sfi'_2$ with distinct tangents:
$T_{\sfi'_1}\cap T_{\sfi'_2}=\{0\}$. Then the degeneration $\sfi'\stackrel{degen}{\to}(\sx^{\!\!\!\!1}\cup\sfi)$ consists of two:
$\sfi'_1\to\sx^{\!\!\!\!1}$ and $\sfi'_2\to\sfi$. Thus the factorization is of the form $(\sx^{\!\!\!\!1}\cup\sx^{\!\!\!\!2})+\sy\ra(\sx^{\!\!\!\!1}\cup\sfi'_2)\ra(\sx^{\!\!\!\!1}\cup\sfi)$.
Finally, the primitivity of $\sx^{\!\!\!\!2}+\sy\ra\sfi$ forces: $\sfi'_2=\sfi$.
\epr

\subsection{Bounds on invariants}%$\empty$\\
\subsubsection{Semi-continuity principle}\label{SecSemiContinuityPrinciple}$\empty$\\
This principle allows to reduce some general questions to the collisions of more restricted types.
\bprop\label{ThmSemiContinuityPrinciple}(cf. proposition \ref{ThmIntroSemicontinuity})
Let $inv$ be an invariant of the singularity type, upper semi-continuous (i.e. non-increasing under the deformations).
\\
\parbox{9cm}
{$\bullet$ Let $S_y\ra S'_y$ be a degeneration and $S_x+S'_y\ra S'_f$ a primitive collision. Then there exists
a primitive collision $S_x+S_y\ra S_f$ and a degeneration, such that the diagram commutes. In particular, $inv(S'_f)\ge inv(S_f)$.
}
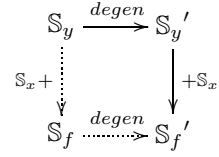
{\xy \POS-(0,-7)\xymatrix{
\sy\ar[r]^{degen} \ar@{.>}[d]_{\sx+}& \sy' \ar[d]^{+\sx}\\
\sfi\ar@{.>}[r]^{degen}&\sfi'
}\endxy}
\li Under the assumptions above, suppose for every primitive collision $S_x+S'_y\ra S'_f$ the bound $inv(S'_f)\le a$ is satisfied.
Then for any primitive collision $S_x+S_y\ra S_f$ one has: $inv(S_f)\le a$.
\eprop
\bpr Note that the degree of curves is assumed to be high. Therefore no pathologies occur, in particular
both $\Si_{\sx\sy}$ and $\Si_\sfi$ are irreducible.

The proof is almost immediate (being just a set theory). The first statement is true because
\mbox{$\ltSi_{\sx\sy}\supset\ltSi_{\sx\sy'}$} causes $\ltSi_{\sx\sy}|_{x=y}\supset\ltSi_{\sx\sy'}|_{x=y}$.
 For the second statement: suppose the degeneration $\sy\ra\sy'$  is done by the intersection
$\ltSi_{\sx\sy}\cap Z=\ltSi_{\sx\sy'}$ (as sets). Then the statement follows from the identity:
\\\mbox{$\Big(\ltSi_{\sx\sy}\cap Z\Big)|_{x=y}=\ltSi_{\sx\sy}|_{x=y}\cap Z$}.
\\\epr
A useful consequence of the principle is the possibility to consider only linear sub-strata. Namely,
let $\Si^{(l)}_{\sx\sy}\subset\Si_{\sx\sy}$ be a linear substratum.
Then \mbox{$\ltSi^{(l)}_{\sx\sy}|_{x=y}\subset\ltSi_{\sx\sy}|_{x=y}$} and all the lower bounds for semi-continuous invariants
of $\Si^{(l)}_{\sx\sy}$ are satisfied for  $\Si_{\sx\sy}$.
$\empty$\\\subsubsection{Multiplicity}\label{SecMultiplicityBounds}$\empty$\\
\bprop
For any initial types $\sx,\sy$ there exists a primitive collision $\sx+\sy\ra \sfi$ with the resulting
multiplicity: $mult(\sfi)=max(mult(\sx),mult(\sy))$.
\eprop
\bpr
Use the semi-continuity principle. First degenerate each of $\sx,\sy$ to a uni-branched \Nnd type
(preserving multiplicities). This can always be done as follows.\\
\parbox{8cm}
{Force all the tangents of a given germ
to coincide. If the so obtained germ is not \Nnd with respect to its \ND, kill all the
necessary monomials, preserving the multiplicity.  (This is always possible by standard arguments from
\cite[section III.3]{AGLV}). }
\begin{picture}(0,0)(-10,20)
\mesh{0}{0}{8}{4}{10}{120}{55}
%% Draws the\mesh{x-coordinate of(0,0),y-coordinate of(0,0),number of x-points(0 isn't counted),number of y-points,delta,x-length,y-length}
\put(0,40){\line(1,-1){10}} \put(10,30){\line(2,-1){20}} \put(30,20){\line(3,-1){30}}\put(60,10){\line(4,-1){40}}
\put(0,40){\line(3,-1){120}}\put(0.5,40){\line(3,-1){120}}
 \put(-2,38){$\bullet$} \put(-10,38){$p$}
\end{picture}
\\
If the so-obtained germ is not \sqh remove the necessary monomials, preserving $x^p_1$.
So, we have arrived to the \sqh germs, of the types $\sx':~~x^{p_x}_1+x^{q_x}_2$, and $\sy':~~x^{p_y}_1+x^{q_y}_2$.

Now collide them such that all the tangents coincide (i.e. $l_x=l=l_y$). Immediate application of the
collision algorithm gives that the multiplicity of the resulting type is $max(mult(\sx'),mult(\sy'))$.
Now invoke the semi-continuity principle.\\
\epr
In general the situation is much more complicated, multiplicity can jump significantly. This happens when the
collision line $l$ and all the non-free tangents are distinct. However there is always the following  bound:
\bprop(cf. proposition \ref{ThmIntroBoundOnMultiplicity})
Let the initial types $\sx\sy$ have the multiplicities $m_x,m_y$ and the numbers of free branches $r_x,r_y$
respectively.
If \mbox{$r_x+r_y\ge m_y$}, then for any collision $\sx+\sy\ra\sfi$: $m_\sfi=m_x$.
If $r_x+r_y< m_y$, then for any collision $\sx+\sy\ra\sfi$: $mult(\sfi)\le m_x-r_x+m_y-r_y$.
\eprop
\bpr The proof goes by first degenerating the types to some specific patterns
(preserving the multiplicities and the number of free branches)
and then applying the semi-continuity principle.
\li Degenerate both $\sx$ and $\sy$ to \gNnd types;
\beq
\sx\ra ~x^{m_x}_1+x^{m_x-r_x}_1x^{r_x}_2+x^{N_x}_2,~~~N_x\gg0,~~~~~\sy\ra ~x^{m_y}_1+x^{m_y-r_y}_1x^{r_y}_2+x^{N_y}_2,~~~N_y\gg0
\eeq
\li By the semi-continuity one can assume both of the degenerated germs to be linear,
i.e. we consider the linear substrata
 $\Si^{(l)}_{\sx'}\subset\Si_{\sx'}$ and  $\Si^{(l)}_{\sy'}\subset\Si_{\sy'}$.
 Thus can write the defining conditions of the stratum $\tSi^{(l)}_{\sx\sy}$
(outside the diagonal $x=y$) explicitly:
\beq\label{EqProofOfMultiplicityProposition}\ber
f|_x^{(m_x+k)}\sim (A_{r_x+k+\de^x_k},\underbrace{l_x.l_x}_{m_x-r_x-\de^x_k}),~~k=0,1...N_x,\\
f|_y^{(m_x+k)}\sim (A_{r_y+m_x-m_y+k+\de^y_k},\underbrace{l_y.l_y}_{m_y-r_y-\de^y_k}),~~k=0,1...N_y
\eer\eeq
So, if $r_x+r_y\ge m_y$, the conditions for $k=0$ can be resolved without increasing the multiplicity:
\beq
f|_x^{(m_x)}\sim (A_{r_x+r_y-m_y},\underbrace{l_x.l_x}_{m_x-r_x},\underbrace{l_y.l_y}_{m_y-r_y-\de^y_k})
\eeq
From the equation (\ref{EqProofOfMultiplicityProposition}) it is seen that all further conditions (with $k>0$) do not increase the multiplicity. So the
final multiplicity is $m_x$.

If  $r_x+r_y<m_y$ then necessarily $f|_x^{(m_x)}=0=f|_x^{(m_x+1)}=...=f|_x^{(m_x+m_y-r_x-r_y-1)}$, while the conditions for
$f|_x^{(m_x+m_y-r_x-r_y)}$ can be resolved in the form\\
$f|_x^{(m_x+m_y-r_x-r_y)}\sim (A_{**},\underbrace{l_x.l_x}_{**},\underbrace{l_y.l_y}_{***})$. As previously, it follows that all
the higher order conditions can be resolved also.\\
\epr
Note that this bound is sharp, e.g. it is realized in the collision of
\mbox{$x^{m_x}_1+x^{m_x-r_x}_1x^{r_x}_2+x^{N_x}_2$}
and \mbox{$x^{m_x}_1+x^{m_y-r_y}_1x^{m_x-m_y+r_y}_2+x^{N_y}_2$} (as in the proof), with $N_x,N_y$ big enough. But it is not the best possible,
e.g. when there are distinct non-free tangents, the bound probably could be improved.

$\empty$\\\subsubsection{How $\de$ changes?}\label{SecDeltaConstCollisions}$\empty$\\
We are particularly interested in $\de=const$ collisions.
By the Milnor-Yung formula $\de=\frac{\mu+r-1}{2}$ and the necessary inequality $\mu_{\sfi}\ge\mu_{\sx}+\mu_\sy+1$
we get immediate
\bprop
%\li Suppose there exists a $\de=const$ collision $\sx+\sy\ra\sfi'$. Then there exists a primitive $\de=const$
%collision that factors the original: $\sx+\sy\ra\sfi\stackrel{degeneration}{\to}\sfi'$
Let $r_x,r_y$ be the (total) number of branches of $\sx\sy$. For a $\de=const$ collision:\\
\mbox{$r_{\sfi}=r_x+r_y-(\mu_\sfi+1-\mu_x-\mu_y)$}. In particular, $r_\sfi \le r_{\sx}+r_{\sy}-2$.
\eprop
 Probably the key result for studying the $\de=const$ collisions is the classical Teissier theorem
\bthe\label{ThmTeissier} \cite{Teissier} The flat family of plane curves $(C_t,0)\to(T,0)$ over a normal base $T$ admits simultaneous normalization iff $\de(C_t)=const$
\ethe
Correspondingly, for the $\de=const$ collision a natural idea is to apply the geometric method as in \S\ref{SecGeometricApproach}:
to blowup at one of points and then to trace the collision on the blown-up plane. First we define a generalization of the notion of free branches.
\bed\label{DefPotentiallyFreeBranches}
Let $C=\cup_i C_i$ be the branch decomposition. A subset $\{C_{i_j}\}_{j\in J}$ is called potentially free if after several blowups
the strict transforms $\tC_{i_j}$ intersect at one point and are free.
\eed
\bex\label{ExamplePotentiallyFreeBranches}
Let $\{C_{i_j}\}$ be smooth branches of constant pairwise tangency, i.e. $deg(C_{i_j}C_{i_k})_{j\ne k}=const$ (independent of $j,k$) and
no other branch intersects them, with higher intersection multiplicity.\\
\parbox{9cm}
{Then after several blowups their strict transforms will
intersect at one point and be pairwise transversal, in addition no other branch will be tangent to any of the chosen branches at this point.
}
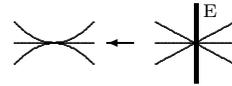
\begin{picture}(0,0)(-30,0)
\curve(-15,0,0,0,15,0) \curve(-15,8,0,0,15,8) \curve(-15,-8,0,0,15,-8)
\put(30,0){\vector(-1,0){10}}
\end{picture}
\begin{picture}(0,0)(-80,0)
\curve(-15,0,0,0,15,0) \curve(-15,8,0,0,15,-8) \curve(-15,-8,0,0,15,8)
\put(0,-15){\line(0,1){30}}\put(0.5,-15){\line(0,1){30}}\put(1,-15){\line(0,1){30}} \put(3,10){\tinyT E}
\end{picture}
\eex
\hspace{-0.6cm}\parbox{10cm}
{A potentially free subset of branches has easy characterization by the resolution tree $\Gamma_\mS$ of the singularity. The
tree contains a subtree as in the picture, where the numbers are the intersection multiplicities with the exceptional divisor.
Denote this subtree by $\Gamma_n$ and its root by $v_n$.
}
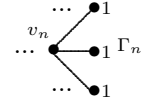
\begin{picture}(0,0)(-30,0)
\put(-14,-1){$...$}  \put(0,15){$...$}   \put(0,-15){$...$}   \put(-9,6){$\tinyA v_n$}    \put(25,0){$\tinyA \Gamma_n$}
\put(-1.5,-2){$\bullet$}\curve(0,0,15,15)\curve(0,0,15,0)\curve(0,0,15,-15)
\put(14,14){$\bullet \tinyM 1$}  \put(14,-2.5){$\bullet\tinyM  1$}  \put(14,-17.5){$\bullet\tinyM  1$}
\end{picture}
\bthe $\bullet$ Assume $m_\sx\ge m_\sy$.
If $\sx$ contains a subset of $m_\sy$ potentially free branches then there exists a collision $\sx+\sy\ra\sfi$ with
the resolution tree $\Gamma_\sfi=(\Gamma_\sx\setminus \Gamma_n)\underset{v_n}{\cup}\Gamma_\sy$, obtained by gluing in
the tree of $\sy$ to the vertex $v_n$ (replacing the subtree $\Gamma_n$).
\li In particular, in such a collision $\de_\sfi=\de_\sx+\de_\sy$, $m_\sfi=m_\sx$,
\mbox{$\mu_\sfi=\mu_\sx+\mu_\sy-1+m_y$}, $r_\sfi=r_x+r_y-m_y$
\ethe
\bpr
Blowup till the potentially free branches become smooth and separated (i.e. one step after the example \ref{ExamplePotentiallyFreeBranches}).
Now glue these smooth branches to $\sy$. To see that this is possible consider a generic line section of $\sy$.
 It intersects the curve with local multiplicity $m_\sy$.
\\
\parbox{8cm}
{Thus deforming the line slightly off the point $y$ gives $m_\sy$ points
of simple (transversal) intersection with the curve.
Therefore the collision is done by moving $\sy$ towards the exceptional divisor (cf. the picture). From this the statement
about the resolution tree follows.
}
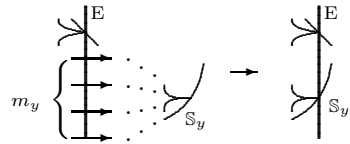
\begin{picture}(0,0)(-30,0)
\curve(0,25,0,-25)\curve(0.5,25,0.5,-25) \put(2,20){\tinyT E}
 \multiput(-5,5)(0,-10){4}{\vector(1,0){15}}  \curve(-10,20,-8,17,0,15)\curve(-10,10,-8,13,0,15)
\curve(-5,20,0,15,5,10)  \put(-28,-13){$\tinyA m_y$ $\Bigg\{$}

\multiput(15,4)(5,-3){3}{.}
\multiput(15,-6)(5,-1){3}{.}\multiput(15,-16)(5,1){3}{.}\multiput(15,-26)(5,3){3}{.}
\curve(30,-5,32,-9,40,-10)\curve(30,-15,32,-11,40,-10)
\curve(30,-18,40,-10,45,3)  \put(38,-20){$\tinyA \sy$}
\put(55,0){\vector(1,0){10}}
\end{picture}
\begin{picture}(0,0)(-115,0)
\curve(0,25,0,-25)\curve(0.5,25,0.5,-25) \put(2,20){\tinyT E}
  \curve(-10,20,-8,17,0,15)\curve(-10,10,-8,13,0,15)
\curve(-5,20,0,15,5,10)

\curve(-10,-5,-8,-9,0,-10)\curve(-10,-15,-8,-11,0,-10)
\curve(-10,-18,0,-10,5,3)  \put(4,-15){$\tinyA \sy$}
\end{picture}
\\
The second statement now follows immediately from the formula $\de=\sum\frac{m_i(m_i+1)}{2}$
 (the summation is over the vertices of the resolution tree, $m_i$ are the multiplicities of the strict transforms) and the formula
 $\mu=2\de-r+1$.
\epr
\subsection{Examples}
$\empty$\\\subsubsection{ADE+ADE$\to$ADE}\label{SecExampleADE}$\empty$\\
By the analysis of Dynkin diagrams and by applying the above algorithm we get the following collisions:
\beq\ber
\xymatrix{
A_k+A_l\ar[r]\ar[rd] & A_{k+l+1}\ar@{~>}[d]
%&\circ\!\!-\!\!\circ\!\!-\!\!\circ\!\!-0-\!\!\circ\!\!-\!\!\circ\!\!-\!\!\circ
\\ & D_{k+l+2}
}~~
\xymatrix{
A_k+A_3\ar[r]\ar[rd]\ar[d] & A_{k+4}\\E_{k+4} & D_{k+4}
}~~
\xymatrix{
A_{k}+A_1\ar[r]\ar[rd] & A_{k+2}\\ & E_{k+2}
}~~
\xymatrix{
A_3+A_2\ar[r]\ar@{~>}[rd]\ar[d] & A_{6}\ar@{~>}[d]\\ D_6\ar@{~>}[r]& E_7
}
\\
\xymatrix{
A_4+A_2\ar[r]\ar@{~>}[rd]\ar[d] & A_{7}\ar@{~>}[d]\\ E_7\ar@{~>}[r]& D_8
}~~
\xymatrix{
D_5+A_k\ar[r]\ar[d] & D_{5+k+1}\\ E_{5+k+1}
}~~
\xymatrix{
E_6+A_1\ar[r]& E_{8}\\D_k+A_l\ar[r]& D_{k+l+1}
}
\eer\eeq
The collisions corresponding to the straight arrows are generic (this can be seen e.g. by codimension or Milnor number).
Wavy arrows indicate the non-generic collision or degeneration. For the types $E_k$, we assume $6\le k\le8$
$\empty$\\\subsubsection{The $D_k$ collisions for some lower cases}$\empty$\\
\beq
D_4+D_4\ra J_{10},~~~
\xymatrix{
D_4+D_5\ar[r]\ar[rd] & X_{1,2},~~\mu=11\\ &J_{2,1},~~\mu=11
}~~
\xymatrix{
D_4+D_6\ar[r]\ar[rd] & X_{1,2},~~\mu=11\\ &J_{2,2},~~\mu=12
}
\eeq
$\empty$\\\subsubsection{Some numerical results}$\empty$\\
The collision of two ordinary multiple points (the minimal lifting):
\\\hspace{-1cm}
\begin{tabular}{|@{}>{$}c<{$}@{}|@{}>{$}c<{$}@{}|@{}>{$}c<{$}@{}|@{}>{$}c<{$}@{}|}\hline
&\sx&\sy&\sfi\\\hline
&x^{p+1}_1+x^{p+1}_2&x^{q+1}_1+x^{q+1}_2&(x^{p-q}_1+x^{p-q}_2)(x^{q+1}_1+x^{2q+2}_2)\\
\mu&p^2&q^2&p^2+q^2+q\\
\de&\frac{p^2+p}{2}&\frac{q^2+q}{2}&\frac{p(p+1)+q(q+1)}{2}\\
\kappa&p^2+p&q^2+q&p^2+p+q^2+q\\\hline
\end{tabular}
\\
Some higher cases, with lifting as in eq. (\ref{EqLiftingSecondAccepted}).
\\
Collision of an ordinary cusp and \omp.
\\
\begin{tabular}{|@{}>{$}c<{$}@{}|@{}>{$}c<{$}@{}|@{}>{$}c<{$}@{}|@{}>{\tinyT$\tinyM}c<{$}@{}|@{}>{\tinyT$\tinyM}c<{$}|@{}>{\tinyT$\tinyM}c<{$}|}\hline
&\sx&\sy&\sfi_{l=l_x}~~p\ge q+2&\sfi_{l=l_x}~~p=q+1&\sfi_{l\ne l_x}\\\hline
&x^{p}_1+x^{p+1}_2&x^{q+1}_1+x^{q+1}_2&(x^{p-1-q}_1+x^{p-q}_2)(x^{q+1}_1+x^{2q+2}_2)&x^{p}_1+x^{2p+1}_2&
(x^{p-q}_1+x^{p-q}_2)(x^{q+1}_1+x^{2q+1}_2)\\
\mu&p^2-p&q^2&p^2-p+(q+1)^2&2p(p-1)&p^2+q^2\\
\de&\frac{p^2-p}{2}&\frac{q^2+q}{2}&\frac{p(p-1)+(q+1)(q+2)}{2}&p(p-1)&\frac{p(p+1)+q(q-1)}{2}\\
\kappa&p^2-1&q^2+q&p^2+q(q+2)&2p^2-p-1&p^2+q^2+q\\\hline
\end{tabular}
\\
Collision of an \omp~ with $\sx=\Big($cusp $\cup$ free branches$\Big):$
\beq\ber
\sx=\prod^r_{i=1} l_i(x^{p}_1+x^{p+1}_2),~~ \mu(\sx)=(p+r-1)^2+p-1,
\\\de(\sx)=\frac{(p+r)(p+r-1)}{2},~~\kappa(\sx)=(p+r)(p+r-1)+p-1
\\
\sy=x^{q+1}_1+x^{q+1}_2,~~q\le p+r,~~\mu(sy)=q^2,~~\de(\sy)=\frac{q^2+q}{2},~~\kappa(\sy)=q^2+q
\eer\eeq
\\\\
\beq\ber
\sfi_{l=l_x}, ~~p\ge q+2,~~ (x^{r}_1+x^{r}_2)(x^{p-q-1}_1+x^{p-q}_2)(x^{q+1}_1+x^{2q+2}_2),
\\
{\mu(\sfi_{l=l_x})=(p+r-1)^2+p-1+(q+1)^2},~~ {\de(\sfi_{l=l_x})=\frac{(p+r)(p+r-1)+(q+1)(q+2)}{2}},
\\
{\kappa(\sfi_{l=l_x})=(p+r)(p+r-1)+p-1+(q+1)^2}
\eer\eeq
\\\\
\beq\ber
\sfi_{l=l_x}, p\le q+1,~~
\sfi_{l=l_x}=(x^{r+p-q-1}_1+x^{r+p-q-1}_2)(x^{q+1-p}_1+x^{2(q+1-p)}_2)(x^{p}_1+x^{2p+1}_2),
\\
{\mu(\sfi_{l=l_x})=(p+r-1)^2+p-1+q^2+q},~~ {\de(\sfi_{l=l_x})=\frac{(p+r)(p+r-1)+q(q+1)}{2}},
\\
{\kappa(\sfi_{l=l_x})=(p+r)(p+r-1)+p-1+q(q+1)}
\eer\eeq
\\\\
\beq\ber
\sfi_{l\ne l_x},~~q\ge r,~~(x^{p+r-q}_1+x^{p+r-q}_2)(x^{q-r+1}_1+x^{2q-2r+1}_2)(x^r_1+x^{2r}_2),\\
\mu(\sfi_{l\ne l_x})=(p+r)^2+q^2-r,~~\de(\sfi_{l\ne l_x})=\frac{(p+r)(p+r+1)+q(q-1)}{2},
\\\kappa(\sfi_{l\ne l_x})=(p+r)^2+p+q^2
\eer\eeq
\\
\beq\ber
\sfi_{l\ne l_x},~~q<r,~~(x^{r-q-1}_1+x^{r-q-1}_2)(x^{p}_1+x^{p+1}_2)(x^{q+1}_1+x^{2(q+1)}_2),\\
\mu(\sfi_{l\ne l_x})=(p+r)^2-p+q(q+1),~~\de(\sfi_{l\ne l_x})=\frac{(p+r)(p+r-1)+q(q+1)+2r}{2},
\\\kappa(\sfi_{l\ne l_x})=(p+r)^2+r+q(q+1)
\eer\eeq

\end{document}